\documentclass{amsart}

\usepackage[english]{babel}

\usepackage[letterpaper,top=3cm,bottom=3cm,left=3cm,right=3cm,marginparwidth=1.75cm]{geometry}

\usepackage{amsmath}
\usepackage{graphicx}
\usepackage{tikz-cd}
\usepackage{amsfonts}
\usepackage{amssymb}
\usepackage{amsthm}
\usepackage[colorlinks=true, allcolors=violet]{hyperref}
\usepackage{cleveref}

\newtheorem{Thm}{Theorem}[section]

\newtheorem{Lm}[Thm]{Lemma}
\newtheorem{Prop}[Thm]{Proposition}

\newtheorem{LetterThm}{Theorem}

\theoremstyle{definition}
\newtheorem{Defn}[Thm]{Definition}

\newtheorem{Ex}[Thm]{Example}
\newtheorem{Rem}[Thm]{Remark} 

\title[Computing the Roots of Bundles arising from Monodromy]{Computing the Roots of Twisting Sheaves over the Projective Line arising from Monodromy Representations}
\author{Diego Yépez}
\thanks{Throughout the majority of this work the author was supported by the Wayne State University Rumble Fellowship.}
\thanks{The author is currently supported by the IDA Postdoctoral Fellowship}
\address{Center for Communications Research, Princeton, NJ 08540}
\email{d.yepez@idaccr.org}
\address{Department of Mathematics, Princeton University, Princeton, NJ 08544}
\email{dy9534@princeton.edu}

\begin{document}

\begin{abstract}
Given a monodromy representation $\rho$ of the projective line minus $m$ points, one can extend the resulting vector bundle with connection map canonically to a vector bundle with logarithmic connection map over all of the projective line. 
Now, since vector bundles split as twisting sheaves over the projective line, the focus of this work regards knowing the exact decomposition; i.e. computing the roots.
Particularly, we compute the roots for all finite-dimensional $\rho$ when $m = 2$ and for all $\rho$ of dimension less than $3$ when $m = 3$. 
\end{abstract}

\date{January 1, 2025}

\maketitle

\tableofcontents

\section{Introduction}
\subsection{Related Problems}

Taking as a point of departure the $m$th punctured Riemann sphere $\mathbb{P}^{1} - \{p_{1},...,p_{m}\}$, it is well-known that a representation of its fundamental group determines a holomorphic connection, see for example \cite{Szamuely2009}. After a choice of logarithm is made, the holomorphic connection can be extended to a unique logarithmic connection over all of $\mathbb{P}^{1}$. Note then, that electing the canonical choice made in \cite{Deligne2006} as a branch of logarithm, a representation of the fundamental group of the $m$th punctured Riemann sphere produces a logarithmic connection on $\mathbb{P}^{1}$. By the celebrated Birkhoff–Grothendieck theorem the underlying vector bundle of the logarithmic connection must decompose as a direct sum of twisting sheaves. The looming question becomes: given a representation of the fundamental group of $\mathbb{P}^{1} - \{p_{1},...,p_{m}\}$ what is the decomposition of the ensuing vector bundle on $\mathbb{P}^{1}$?  

First, observe that in the case of a compact (no points removed) Riemann surface, the question chiefly translates to which vector bundles admit holomorphic connections. This was settled by Weil \cite{Weil1938}: a holomorphic vector bundle on a Riemann surface admits a holomorphic connection if and only if each indecomposable summand has first Chern class 0.

Of course, as the question under scrutiny deals with monodromy representations and logarithmic connections, a natural consideration is to look at the moduli of such objects. While the moduli spaces in question have been studied for sufficiently nice varieties, for example by Simpson in \cite{Simpson1994} and Nitsure in \cite{Nitsure1993}, the results of such papers, albeit extremely important in their own right, are too general to produce an answer to our specific question. However, if we restrict the problem to only consider unitary monodromy representations of $\mathbb{P}^{1} - \{p_{1},...,p_{m}\}$, then by the work of Mehta and Seshadri in \cite{MS1980}, this essentially amounts to the study of the moduli space of parabolic vector bundles on $\mathbb{P}^{1}$ with logarithmic connection map; particularly, Mehta and Seshadri showed that there is a one-to-one correspondence between irreducible unitary monodromy representations and stable parabolic vector bundles with logarithmic
connection map. 

Regarding general facts of the moduli of parabolic vector bundles on $\mathbb{P}^{1}$ with logarithmic
connection map one may look for example at the following papers \cite{Inaba2006, Inaba2007, DP2022, Fer2021}. Of particular interest to us, are the result of \cite{Hu2024}, \cite{Mats2024} and \cite{MT2021}. 
The authors of \cite{Hu2024} generalize some of the results from \cite{Loray2015, Loray2013} and provide the decomposition of any rank $2$ stable parabolic vector bundle of degree $d$ over $\mathbb{P}^{1}$ with logarithmic connection map over $4$ or more marked points. 
In our setting, this means that the answer to our question is known when we remove $4$ or more points and take $\rho$ to be an irreducible unitary representation of dimension $2$. In \cite{Mats2024}, Takafumi finds the decomposition of any rank $3$ stable parabolic vector bundle of degree $-2$ over $\mathbb{P}^{1}$ with logarithmic connection map over $3$ marked points. The degree $-2$ condition translates to the sum of the traces of the principal logarithm of the image of the generators of the fundamental group being equal to $2$. 
The authors of \cite{MT2021} find a bound for all twisting parameters, $\xi_{j}$, in the decomposition of a $\text{rank} \geq 2$ stable parabolic vector bundle with logarithmic connection map with an extra set of conditions; they state that for each $\xi_{j}$, $-m < \xi_{j} < 0$. Again, in our setting, this means that the answer to our question is bounded when we remove $m$ points and take $\rho$ to be an irreducible unitary representation of dimension two or greater that satisfies the extra set of conditions.  

In this paper, we are not necessarily interested in the overall description of certain moduli spaces but simply in the computation of the decomposition of the vector bundles over the projective line arising from a monodromy representation, as a result we take a rather algebraic approach and consider generic representations. Particularly, we compute the decomposition for all finite-dimensional $\rho$ when $m = 2$ and for all representations $\rho$ of dimension less than $3$ when $m = 3$. Lastly, we remark that we use different terminology than the one use in the setting of moduli of parabolic bundles on $\mathbb{P}^{1}$ with logarithmic
connection map.

\subsection{Results}

Let $\textbf{Y}_{(p_{1}, ..., p_{m})} := \mathbb{P}^{1}_{\mathbb{C}} - \{p_{1}, ..., p_{m}\}$. 
The fundamental group of $\textbf{Y}_{(p_{1}, ..., p_{m})}$ is known to be isomorphic to the free group on $m-1$ generators, $\mathbb{Z}^{*(m-1)}$. 
Thus, a representation of the fundamental group of $\textbf{Y}_{(p_{1}, ..., p_{m})}$ is determined by the images of the generators, each of which we call a \textit{local monodromy}. 
A \textit{monodromy representation} of $\textbf{Y}_{(p_{1}, ..., p_{m})}$ is a representation of its fundamental group. 

We denote a \textit{holomorphic connection} on $\textbf{Y}_{(p_{1}, ..., p_{m})}$ with a pair $(\mathcal{V}, \nabla)$ where $\mathcal{V}$ is a locally free sheaf on $\textbf{Y}_{(p_{1}, ..., p_{m})}$, and $\nabla: \mathcal{V} \rightarrow \mathcal{V} \otimes \Omega^{1}_{\textbf{Y}_{(p_{1}, ..., p_{m})}}$ is a morphism of sheaves satisfying the Leibniz Rule.
When we want to refer to $\nabla$ alone we call it the \textit{connection map}.

Let $S$ be a finite number of points of $\mathbb{P}^{1}$ and $\Omega_{\mathbb{P}^{1}}^{1}(S)$ the sheaf of $1$-forms with logarithmic poles along $S$. 
We denote a \textit{connection with logarithmic poles along $S$} on $\mathbb{P}^{1}$ with a pair $(\overline{\mathcal{E}}, \overline{\nabla})$, where $\overline{\mathcal{E}}$ is a locally free sheaf on $\mathbb{P}^{1}$, and $\overline{\nabla}: \overline{\mathcal{E}} \rightarrow \overline{\mathcal{E}} \otimes_{\mathcal{O}} \Omega_{\mathbb{P}^{1}}^{1}(S)$ is a $\mathbb{C}$-linear map satisfying the Leibniz Rule.

There is a equivalence of categories between finite-dimensional monodromy representations, $\rho$, of $\textbf{Y}_{(p_{1}, ..., p_{m})}$, and holomorphic connections, $(\mathcal{V}, \nabla)$, on $\textbf{Y}_{(p_{1}, ..., p_{m})}$. 
Furthermore, we can extend a holomorphic connection $(\mathcal{V}, \nabla)$ on $\textbf{Y}_{(p_{1}, ..., p_{m})}$, to a connection $(\overline{\mathcal{V}}, \overline{\nabla})$ on $\mathbb{P}^{1}$ with logarithmic poles along $\{p_{1}, ... ,p_{m} \}$ satisfying $(\overline{\mathcal{V}}, \overline{\nabla})|_{\textbf{Y}_{(p_{1}, ..., p_{m})}} \cong (\mathcal{V}, \nabla)$ such that the residue at each pole is the principal logarithm of the local monodromies. 
If $(\mathcal{V},\nabla)$ is the homolomorphic connection associated to $\rho$ then we refer to $(\overline{\mathcal{V}}, \overline{\nabla})$ as the \textit{associated extended logarithmic connection} and denote it as $(\mathcal{V}_{\text{Log}(\rho)}, \nabla_{\text{Log}(\rho)})$.

Thus, the question at hand is, given a monodromy representation of $\textbf{Y}_{(p_{1}, ..., p_{m})}$, extending the resulting holomorphic connection canonically to a connection with logarithmic poles on $\mathbb{P}^{1}$, and applying the Birkhoff–Grothendieck theorem, what is the explicit decomposition of the ensuing vector bundle as a direct sum of twisting sheaves? 
In this context, we refer to the twisting parameters of the twisting sheaves as the \textit{roots}.

In the case of $\textbf{Y}_{(0,\infty)}$ we can fully answer our question. 
When the monodromy representation is a character, the root is the first Chern class.
Given a locally free sheaf $\mathcal{E}$ on $\mathbb{P}^{1}$, denote its first Chern class by $c_{1}(\mathcal{E})$. 
\begin{LetterThm}[\Cref{RootGen}] \label{Red-Z}
    Let $\rho: \mathbb{Z} \rightarrow \text{GL}_{n}(\mathbb{C})$ be an $n$-dimensional monodromy representation of $\textbf{Y}_{(0,\infty)}$ with $(\mathcal{V}_{Log(\rho)}, \nabla_{Log(\rho)})$ the associated extended logarithmic connection. Allow $\zeta = c_{1}(\mathcal{V}_{\text{Log}(\rho)})$.
    \begin{enumerate}
        \item If $\rho$ is an indecomposable representation determined by the Jordan block $M$ with eigenvalue $\lambda = re^{2 \pi iq}$ where $0 \leq q < 1$, then 
        \begin{eqnarray} 
        \nonumber
        \mathcal{V}_{Log(\rho)} \cong 
        \begin{cases} 
        \mathcal{O}(0)^{\oplus n} & \text{when } q=0 \\
        \mathcal{O}(-1)^{\oplus n} & \text{when } q \not= 0.
        \end{cases}
        \end{eqnarray}
        \item If $\rho \cong \bigoplus_{i = 1}^{u} \rho_{i}$ where each $\rho_{i}$ is indecomposable, then 
        \begin{eqnarray}
        \nonumber
        \mathcal{V}_{Log(\rho)} \cong \mathcal{O}(-1)^{\oplus -\zeta} \oplus \mathcal{O}(0)^{\oplus n + \zeta}.
        \end{eqnarray}
        \end{enumerate}
\end{LetterThm}

\noindent Observe that the theorem above closes the case on $\textbf{Y}_{(0,\infty)}$. As \Cref{Ohtsuki} shows, $\text{c}_{1}(\mathcal{V}_{Log(\rho)})$ is determined by $\rho$. Hence, \Cref{Red-Z}$(2)$ reduces finding the roots to computing the first Chern class. 

With regard to $\textbf{Y}_{(0,1,\infty)}$, we restrict ourselves to only consider monodromy representations of dimension less than $3$. 
When the monodromy representation is a character, the root is the first Chern class. 
When $\rho$ is a reducible $2$-dimensional monodromy representation of $\textbf{Y}_{(0,1,\infty)}$ we obtain the following theorem: 

\begin{LetterThm}[\Cref{NeedDer}] \label{Red-two}
    Suppose that $0 \rightarrow \rho' \rightarrow \rho \rightarrow \rho'' \rightarrow 0$ is a short exact sequence of monodromy representations of $\textbf{Y}_{(0, 1, \infty)}$ with $\rho$ two-dimensional and $\rho', \rho''$ one-dimensional. 
    Let $0 \rightarrow \mathcal{O}(\xi') \rightarrow \mathcal{V}_{Log(\rho)} \rightarrow \mathcal{O}(\xi'') \rightarrow 0$ be the short exact sequence of the associated extended vector bundles with logarithmic connection maps.
    \begin{enumerate}
        \item If $\xi',\xi'' \not= -2,0$ respectively, then $0 \rightarrow \mathcal{O}(\xi') \rightarrow \mathcal{V}_{Log(\rho)} \rightarrow \mathcal{O}(\xi'') \rightarrow 0$ splits.
        \item If $\xi', \xi'' = -2,0$ respectively, then $\mathcal{V}_{Log(\rho)} \cong \mathcal{O}(-1) \oplus \mathcal{O}(-1)$ or $\mathcal{V}_{Log(\rho)} \cong \mathcal{O}(0) \oplus \mathcal{O}(-2)$.
    \end{enumerate}  
\end{LetterThm}

In contrast with $\textbf{Y}_{(0, \infty)}$, the fundamental group of $\textbf{Y}_{(0,1,\infty)}$ admits irreducible representations of dimension greater than $1$.
Now, when $\rho$ is an irreducible representation of $\mathbb{Z}*\mathbb{Z}$ the methods we have been using are not sufficient to find the roots as explained at the beginning of $\S 5$, so we are in need of a new tool. Observing the use of the Serre derivative in \cite[\S $3$ and \S $4$]{CHMY2018}, we set out to construct a Serre-like derivative operator. 

Take the character representation $\chi: \pi_{1}(\mathbf{Y}_{(0,1,\infty)}, y) \rightarrow \mathbb{C}^{*}$ defined by mapping both generators to $-1$, the resulting associated extended logarithmic connection is $(\mathcal{V}_{Log(\chi)}, \nabla_{Log(\chi)})$ over $\mathbb{P}^{1}$ with $\mathcal{V}_{Log(\chi)} \cong \mathcal{O}(-1)$. 
Next, we dualize, and obtain $(\mathcal{O}(1), \nabla)$ with 
\begin{eqnarray}
    \nabla: \mathcal{O}(1) \rightarrow \mathcal{O}(1) \otimes_{\mathcal{O}} \Omega_{\mathbb{P}^{1}}^{1}([0]+[1]+[\infty])
\end{eqnarray}
\noindent where  
\begin{eqnarray}
    \Omega_{\mathbb{P}^{1}}^{1}([0]+[1]+[\infty]) \cong \mathcal{O}(1)
\end{eqnarray}

\noindent and so
\begin{eqnarray}
    \nabla: \mathcal{O}(1) \rightarrow \mathcal{O}(2).
\end{eqnarray}

Further, we take $\nabla^{\otimes 2}$ to be the connection map on $\mathcal{O}(1) \otimes \mathcal{O}(1)$ such that for all sections $s,t \in O(1)$ the following is satisfied: $(\nabla \otimes \nabla)(s \otimes t) = \nabla(s) \otimes t + s \otimes \nabla(t)$. 
Observe that 
\begin{eqnarray}
    \nabla^{\otimes 2}: \mathcal{O}(2) \rightarrow \mathcal{O}(2) \otimes \Omega_{\mathbb{P}^{1}}^{1}([0]+[1]+[\infty]),
\end{eqnarray}
 
\noindent and so, inductively,  we define $\nabla_{\xi} := \nabla^{\otimes \xi} : \mathcal{O}(\xi) \rightarrow \mathcal{O}(\xi) \otimes_{\mathcal{O}} \Omega^{1}([0]+[1]+[\infty])$. 
Hence, using the fact that $ \Omega_{\mathbb{P}^{1}}^{1}([0]+[1]+[\infty]) \cong \mathcal{O}(1)$, then
\begin{eqnarray}
    \nabla_{\xi}: \mathcal{O}(\xi) \rightarrow \mathcal{O}(\xi+1).
\end{eqnarray}

\noindent We refer to $\nabla_{\xi}$ as the \textit{auxiliary connection map of weight} $\xi$ \textit{arising from} $\textbf{Y}_{(0,1,\infty)}$.
 
Thus, keeping in mind that $\Omega_{\mathbb{P}^{1}}^{1}([0]+[1]+[\infty]) \cong \mathcal{O}(1)$ we now have the extra tool 
\begin{eqnarray}
    \nabla_{Log(\rho)} \otimes \nabla: \mathcal{V}_{Log(\rho)} \otimes_{\mathcal{O}} \mathcal{O}(1) \rightarrow \mathcal{V}_{Log(\rho)} \otimes_{\mathcal{O}} \mathcal{O}(2) 
\end{eqnarray}

\noindent and more generally, considering an auxiliary connection map of arbitrary weight 
\begin{eqnarray}
    \nabla_{Log(\rho)} \otimes \nabla_{\xi}: \mathcal{V}_{Log(\rho)} \otimes_{\mathcal{O}} \mathcal{O}(\xi) \rightarrow \mathcal{V}_{Log(\rho)} \otimes_{\mathcal{O}} \mathcal{O}(\xi+1). 
\end{eqnarray}
Let $\mathcal{N}_{\xi}(\rho) := H^{0}(\mathbb{P}^{1}, \mathcal{V}_{\text{Log}(\rho)} \otimes_{\mathcal{O}} \mathcal{O}(\xi))$, and furthermore define $\mathcal{N}(\rho) := \bigoplus_{\xi \in \mathbb{Z}} \mathcal{N}_{\xi}(\rho)$. 
We refer to $\mathcal{N}(\rho)$ as the \textit{twisted module of $\rho$}.

Moreover, $\mathcal{N}(\rho)$ is a $\mathbb{Z}$-graded module of global sections over the ring  $\mathbb{C}[x,y]$ where of course $\mathbb{P}^{1} = \text{Proj } \mathbb{C}[x,y]$. 
On global sections, the new connection maps $\nabla_{Log(\rho)} \otimes \nabla_{\xi}$ give us a graded derivation of degree one. 
Indeed, let us denote the connection maps $\nabla_{Log(\rho)} \otimes \nabla_{\xi}$ as $D_{\xi}$ when applied to global sections, so that $D_{\xi}: \mathcal{N}_{\xi}(\rho) \rightarrow \mathcal{N}_{\xi+1}(\rho) $. Observe that each $D_{\xi}$ satisfies the Leibniz rule as each $\nabla_{Log(\rho)} \otimes \nabla_{\xi}$ is a connection map and so 
\begin{eqnarray}
    D := \bigoplus_{\xi \in \mathbb{Z}} D_{\xi}
\end{eqnarray}

\noindent is a graded derivation of degree one on $\mathcal{N}(\rho)$.
We refer to $D$ as the \textit{monodromy derivative of} $\rho$.

With the use of the monodromy derivative of $\rho$ we are able to prove the following theorem:
\begin{LetterThm}[\Cref{even-odd}]
    Suppose that $\rho: \mathbb{Z}*\mathbb{Z} \rightarrow \text{GL}_{2}(\mathbb{C})$ is an irreducible two-dimensional monodromy representation of $\textbf{Y}_{(0,1,\infty)}$ and let $(\mathcal{V}_{\text{Log}(\rho)}, \nabla_{\text{Log}(\rho)})$ denote the associated extended logarithmic connection. 
    Allow $\zeta = c_{1}(\mathcal{V}_{\text{Log}(\rho)})$.
    Then 
    \begin{eqnarray}
        \nonumber
        \mathcal{V}_{\text{Log}(\rho)} \cong
        \begin{cases}
        \mathcal{O}(\frac{\zeta}{2})^{\oplus 2} & \text{when } \zeta \text{ is even} \\
        \mathcal{O}(\frac{\zeta - 1}{2}) \oplus \mathcal{O}(\frac{\zeta + 1}{2}) & \text{when } \zeta \text{ is odd}. 
    \end{cases}
    \end{eqnarray} 
\end{LetterThm}

Witness that the theorem above does not make any assumptions regarding $\rho$ aside from the fact that it is an irreducible representation of dimension $2$.
Thus, this theorem along with \Cref{Red-two} close the case on two-dimensional monodromy representations of $\textbf{Y}_{(0,1,\infty)}$.
$\vspace{3mm}$

$\textbf{Conventions and Assumptions}$. 
\begin{itemize}
    \item All representations $\rho$ are assumed to be over $\mathbb{C}$.
    \item Whenever a branch of the logarithm is needed we will always take the principal branch even when not explicitly stated i.e. for all $re^{2\pi i q} \in \mathbb{C}, 0 \leq q < 1$.
    \item We denote the generators of the free group on two generators with $\gamma_{0}$ and $\gamma_{1}$.
\end{itemize}
$\vspace{1mm}$

$\textbf{Acknowledgments}$. I am greatly indebted to Luca Candelori for all the insightful discussions. I'd like to thank Nick Rekuski for both a careful reading of a previous version of this work and for all the helpful comments. I'd also like to thank Nick Katz and the Mathematics Department at Princeton University for their hospitality.

\section{Preliminaries}

The purpose of this section is to formally provide some background, introduce notation, and state a special case of the main result from \cite{Ohtsuki1982}, which is a fundamental tool in this work. 

Let $\textbf{Y}_{(p_{1}, ..., p_{m})} := \mathbb{P}^{1} - \{p_{1}, ..., p_{m}\}$. 
The fundamental group of $\textbf{Y}_{(p_{1}, ..., p_{m})}$ is known to be isomorphic to the free group on $m-1$ generators, $\mathbb{Z}^{*(m-1)}$. 
Thus, a representation of the fundamental group of $\textbf{Y}_{(p_{1}, ..., p_{m})}$ is determined by the images of the generators, each of which we call a \textit{local monodromy}. 

\begin{Defn}
    A \textit{monodromy representation} of $\textbf{Y}_{(p_{1}, ..., p_{m})}$ is a representation of its fundamental group. 
\end{Defn}

\begin{Defn}
    A \textit{holomorphic connection} on $\textbf{Y}_{(p_{1}, ..., p_{m})}$ is a pair $(\mathcal{V}, \nabla)$ where $\mathcal{V}$ is a locally free sheaf on $\textbf{Y}_{(p_{1}, ..., p_{m})}$, and $\nabla: \mathcal{V} \rightarrow \mathcal{V} \otimes \Omega^{1}_{\textbf{Y}_{(p_{1}, ..., p_{m})}}$ is a morphism of sheaves satisfying the Leibniz Rule.
    When we want to refer to $\nabla$ alone we call it the \textit{connection map}.
    A morphism of connections $(\mathcal{V}, \nabla) \mapsto (\mathcal{V}', \nabla')$ is a morphism of $\mathcal{O}$-modules $\phi: \mathcal{V} \rightarrow \mathcal{V}'$ making the diagram 
    \begin{center}
    \begin{tikzcd}
    \mathcal{V} \arrow[r, "\nabla"] \arrow[d, "\phi"]
    & \mathcal{V} \otimes \Omega^{1}_{\textbf{Y}_{(p_{1}, ..., p_{m})}} \arrow[d, "\phi \otimes \text{id}"] \\
    \mathcal{V}' \arrow[r, "\nabla'"]
    & \mathcal{V} \otimes \Omega^{1}_{\textbf{Y}_{(p_{1}, ..., p_{m})}} 
    \end{tikzcd}
    \end{center}
    \noindent commute.
\end{Defn}

\begin{Thm} \label{Free-sheaf}
    Every locally free sheaf on $\textbf{Y}_{(p_{1}, ..., p_{m})}$ is free.
\end{Thm}

\proof See \cite[Theorem $30.4$]{Forster2012}.
$\Box$

\begin{Lm} \label{Equiv-cat}
    The category of finite-dimensional monodromy representations of $\textbf{Y}_{(p_{1}, ..., p_{m})}$ is equivalent to the category of holomorphic connections on $\textbf{Y}_{(p_{1}, ..., p_{m})}$.
\end{Lm}

\proof Combine \cite[Corollary $2.6.2$]{Szamuely2009} with \cite[Proposition $2.7.5$]{Szamuely2009}.
$\Box$

\begin{Defn} \label{Associated-Con}
    By \Cref{Equiv-cat}, given a monodromy representation of $\textbf{Y}_{(p_{1}, ..., p_{m})}$, $\rho$, we obtain a holomorphic connection on $\textbf{Y}_{(p_{1}, ..., p_{m})}$. We call such a connection the $\textit{associated complex connection}$ and denote it by $(\mathcal{Q}_{\rho}, \nabla_{\rho})$.
\end{Defn}

\begin{Defn}
    Let $S$ be a finite number of points of $\mathbb{P}^{1}$ and $\Omega_{\mathbb{P}^{1}}^{1}(S)$ the sheaf of $1$-forms with logarithmic poles along $S$. 
    A \textit{connection with logarithmic poles along S} on $\mathbb{P}^{1}$ is a pair $(\overline{\mathcal{E}}, \overline{\nabla})$, where $\overline{\mathcal{E}}$ is a locally free sheaf on $\mathbb{P}^{1}$, and $\overline{\nabla}: \overline{\mathcal{E}} \rightarrow \overline{\mathcal{E}} \otimes_{\mathcal{O}} \Omega^{1}(S)$ is a $\mathbb{C}$-linear map satisfying the Leibniz Rule.
\end{Defn}

\begin{Lm} \label{Extension}
    Given a holomorphic connection $(\mathcal{V}, \nabla)$ on $\textbf{Y}_{(p_{1}, ..., p_{m})}$, there is a unique connection $(\overline{\mathcal{V}}, \overline{\nabla})$ on $\mathbb{P}^{1}$ with logarithmic poles along $\{p_{1}, ... ,p_{m} \}$ satisfying $(\overline{\mathcal{V}}, \overline{\nabla})|_{\textbf{Y}_{(p_{1}, ..., p_{m})}} \cong (\mathcal{V}, \nabla)$ such that the residue at each pole is the principal logarithm of the local monodromies.
    Moreover, the extension $(\overline{\mathcal{V}}, \overline{\nabla})$ of $(\mathcal{V}, \nabla)$ is functorial and exact in $(\mathcal{V}, \nabla)$.
\end{Lm}

\proof See \cite[Proposition $5.4$]{Deligne2006}.
$\Box$

\begin{Defn} \label{extended-log-definition}
    As a result of applying \Cref{Extension} to \Cref{Equiv-cat} we can associate to each monodromy representation $\rho$ of $\textbf{Y}_{(p_{1}, ..., p_{m})}$ a unique logarithmic connection over $\mathbb{P}^{1}$. 
    We refer to such connections as the $\textit{associated extended logarithmic connection}$ and denote it as $(\mathcal{V}_{\text{Log}(\rho)}, \nabla_{\text{Log}(\rho)})$.
\end{Defn}

\begin{Defn} \label{roots}
    In the setting of \Cref{extended-log-definition}, as $\mathcal{V}_{\text{Log}(\rho)} \cong \bigoplus_{j = 1}^{\text{dim} (\rho)} \mathcal{O}(\xi_{j})$ we call the twisting parameters $\xi_{j}$ the \textit{roots}.
\end{Defn}

\begin{Lm} \label{first-chern-class}
    The first Chern class of $\mathcal{V} \in \textnormal{Vect}(\mathbb{P}^{1})$, denoted $c_{1}(\mathcal{V})$, is the sum of the twisting parameters. 
\end{Lm}

\proof As every vector bundle over $\mathbb{P}^{1}$ splits as a direct sum of twisting sheaves, it follows that $\mathcal{V} \cong \bigoplus\limits_{j = 1}^{r} \mathcal{O}(\xi_{j})$, and so det$(\mathcal{V}) = \text{det}(\bigoplus\limits_{j = 1}^{r} \mathcal{O}(\xi_{j})) = \bigotimes\limits_{j = 1}^{r} \mathcal{O}(\xi_{j}) = \mathcal{O}(\sum\limits_{j = 1}^{r} \xi_{j})$. Thus, $\text{c}_{1}(\mathcal{V}) = \text{deg}(\text{det}(\mathcal{V})) = \text{deg}(\mathcal{O}(\sum\limits_{j = 1}^{r} \xi_{j})) = \sum\limits_{j = 1}^{r} \xi_{j}$.
$\Box$

\begin{Thm}
    Let $(\mathcal{V}, \nabla)$ be a connection on $\mathbb{P}^{1}$ with logarithmic poles along $\{p_{1}, ..., p_{m}\}$, then the following holds:
    \begin{eqnarray}
    \label{Ohtsuki}
    \text{c}_{1}(\mathcal{V}) = - \sum\limits_{j=1}^{m} \text{Tr}(\text{Res}_{p_{j}}\nabla)
\end{eqnarray}

\noindent where $\text{Res}_{p_{j}}\nabla$ is the residue of $\nabla$ at the pole $p_{j}$. 
\end{Thm}

\proof See \cite{Ohtsuki1982}.

\begin{Rem}
    Observe that \Cref{Ohtsuki} provides for us an explicit formula for computing the first Chern class of any associated extended logarithmic connection. 
    Moreover, observe that any higher Chern class vanishes as we are over $\mathbb{P}^{1}$.
\end{Rem}

\begin{Defn}
   We refer to Equation $(\ref{Ohtsuki})$ as $\textit{Ohtsuki's formula}$.  
\end{Defn} 

\begin{Rem}
    We note that \Cref{Free-sheaf}, \Cref{Extension} and \Cref{Ohtsuki} are special cases of the much more general statements proved in their respective sources.
\end{Rem}

\section{Punctured Complex Line}

The goal of this section is to compute the roots arising from any monodromy representation of $\textbf{Y}_{(0, \infty)}$. 

\begin{Prop} \label{RootC}
    Let $\chi: \mathbb{Z} \rightarrow \mathbb{C}^{*}$ be a monodromy representation of $\textbf{Y}_{(0,\infty)}$ given by $1 \mapsto re^{2\pi iq}$ with $0 \leq q < 1$, and $(\mathcal{V}_{Log(\chi)}, \nabla_{Log(\chi)})$ the associated extended logarithmic connection. Then
    \begin{eqnarray} \label{Char-Root}
    \nonumber
    \mathcal{V}_{Log(\chi)} \cong 
    \begin{cases} 
      \mathcal{O}(0) & \text{when } q=0 \\
      \mathcal{O}(-1) & \text{when } q \not= 0.
    \end{cases}
    \end{eqnarray}
\end{Prop}

\proof Observe that the connection map has two poles, one at $0$ and one at $\infty$ with respective residues $\text{log} (re^{2 \pi iq})$ and $\text{log} (\frac{1}{r}e^{-2 \pi iq})$.

By Ohtsuki's formula we have that $c_{1}(\mathcal{V}_{Log(\chi)}) = -\text{log} (re^{2 \pi iq}) -\text{log} (\frac{1}{r}e^{-2 \pi iq})$, choosing the branch of log defined by $0 \leq q < 1$ we see that if $q = 0$, then $c_{1}(\mathcal{V}_{Log(\chi)}) = -\text{log}(r) - 0 - \text{log}(r^{-1}) + 0$. 
However, since $r \in \mathbb{R}$ then $\text{log}(r^{-1}) = -\text{log}(r)$ and so $c_{1}(\mathcal{V}_{Log(\chi)}) = -\text{log}(r) + \text{log}(r) = 0$.
If $q \not= 0$, as $-q \not\in [0,1)$, we must use $-q+1$, thus, $c_{1}(\mathcal{V}_{Log(\chi)}) = -\text{log}(r) -q +\text{log}(r) +q - 1 = -1$. 

Hence, 
\begin{eqnarray} 
\nonumber
c_{1}(\mathcal{V}_{Log(\chi)}) = 
 \begin{cases} 
      0 & \text{when } q=0 \\
      -1 & \text{when } q \not= 0
 \end{cases}
\end{eqnarray}

\noindent and the result follows as a consequence of \Cref{first-chern-class}. 
$\Box$

\begin{Lm} \label{ChernR}
     Suppose that $\rho$ is a monodromy representation of $\textbf{Y}_{(0,\infty)}$ defined by $1 \mapsto M$ with Jordan decomposition by blocks $M_{j}$ with eigenvalue $\lambda_{j} = re^{2\pi iq_{j}}$ where $0 \leq q_{j} < 1$. The first Chern class of the associated extended logarithmic connection $(\mathcal{V}_{Log(\rho)}, \nabla_{Log(\rho)})$ is 
     $$c_{1}(\mathcal{V}_{Log(\rho)}) = \sum\limits_{j =1}^{\alpha} \text{dim}(M_{j}) c_{1}(\mathcal{V}_{\text{Log}(\chi_{q_{j}})})$$ 
     
     \noindent where $\alpha$ is the number of Jordan blocks that make up $M$. 
\end{Lm}

\proof Recall the fact that if $S \in \text{GL}_{n}(\mathbb{C})$ and $T \in \text{Mat}(n \times n, \mathbb{C})$ then
\begin{eqnarray} \label{Jf}
    \text{exp}(S^{-1}TS) = S^{-1}\text{exp}(T)S.
\end{eqnarray}
Observe that as a consequence we may suppose that $M$ is in Jordan form. 
Now, take an arbitrary Jordan block $M_{j}$ with eigenvalue $\lambda_{j} = re^{2\pi iq_{j}}$. 
Note that by assumption $q_{j} \in [0,1)$. Then since $\text{log}(M_{j}) = \text{log}(\lambda_{j})\mathbb{I} + N$ where $N$ is a nilpotent matrix, we have $\text{Tr}(\text{log}(M_{j})) = \sum\limits_{1}^{\text{dim}(M_{j})} \text{log}(\lambda_{j})$. 
Next, since we are only interested in the trace, we can choose a basis where $M_{j}^{-1}$ is a Jordan matrix to easily obtain the trace of $\text{log}(M_{j}^{-1})$.
It is clear that the eigenvalue of $M_{j}^{-1}$ is $\lambda^{-1}_{j}$, and so 
\begin{eqnarray}
    \nonumber 
    \text{Tr}(\text{log}(M_{j}^{-1})) = \sum\limits_{1}^{\text{dim}(M_{j})} \text{log}(\lambda^{-1}_{j}).
\end{eqnarray}

Observe that $M_{j}$ determines a subrepresentation $\rho_{j}$, which by \Cref{Equiv-cat} and \Cref{Extension} we have $(\mathcal{V}_{\text{Log}(\rho_{j})}, \nabla_{\text{Log}(\rho_{j})})$. 
Thus, by Ohtsuki's formula we may write
\begin{eqnarray}
\nonumber c_{1}(\mathcal{V}_{Log(\rho_{j})}) &&= -\text{Tr}(\text{log}(M_{j})) - \text{Tr}(\text{log}(M_{j}^{-1})) \\
\nonumber && =  \sum\limits_{1}^{\text{dim}(M_{j})} -\text{log}(\lambda_{j}) - \sum\limits_{1}^{\text{dim}(M_{j})} \text{log}(\lambda^{-1}_{j})\\
\nonumber && = \sum\limits_{1}^{\text{dim}(M_{j})} -\text{log}(\lambda_{j}) - \text{log}(\lambda^{-1}_{j}) \\
\nonumber && = \sum\limits_{1}^{\text{dim}(M_{j})} c_{1}(\mathcal{V}_{\text{Log}(\chi_{q_{j}})}) \\
\nonumber && = \text{dim}(M_{j}) c_{1}(\mathcal{V}_{\text{Log}(\chi_{q_{j}})}).
\end{eqnarray}

\noindent If $M$ is an indecomposable Jordan block then we are done. If not, then it is clear that $c_{1}(\mathcal{V}_{Log(\rho)}) = \sum\limits_{j =1}^{\alpha} \text{dim}(M_{j}) c_{1}(\mathcal{V}_{\text{Log}(\chi_{q_{j}})})$ where $\alpha$ is the number of Jordan blocks that make up $M$.
$\Box$

\begin{Thm} \label{RootGen}
    Let $\rho: \mathbb{Z} \rightarrow \text{GL}_{n}(\mathbb{C})$ be an $n$-dimensional monodromy representation of $\textbf{Y}_{(0,\infty)}$ with $(\mathcal{V}_{Log(\rho)}, \nabla_{Log(\rho)})$ the associated extended logarithmic connection. Allow $\zeta = c_{1}(\mathcal{V}_{\text{Log}(\rho)})$.
    \begin{enumerate}
        \item If $\rho$ is an indecomposable representation determined by the Jordan block $M$ with eigenvalue $\lambda = re^{2 \pi iq}$ where $0 \leq q < 1$, then 
        \begin{eqnarray} 
        \nonumber
        \mathcal{V}_{Log(\rho)} \cong 
        \begin{cases} 
        \mathcal{O}(0)^{\oplus n} & \text{when } q=0 \\
        \mathcal{O}(-1)^{\oplus n} & \text{when } q \not= 0.
        \end{cases}
        \end{eqnarray}
        \item If $\rho \cong \bigoplus_{i = 1}^{u} \rho_{i}$ where each $\rho_{i}$ is indecomposable, then 
        \begin{eqnarray}
        \nonumber
        \mathcal{V}_{Log(\rho)} \cong \mathcal{O}(-1)^{\oplus -\zeta} \oplus \mathcal{O}(0)^{\oplus n + \zeta}.
        \end{eqnarray}
        \end{enumerate}
\end{Thm}

\proof For part $(1)$, we proceed by induction on the dimension of $\rho$.

$\underline{\text{Base Case}}$. Beginning with a character $\chi: \mathbb{Z} \rightarrow \mathbb{C}^{*}$ determined by $1 \mapsto re^{2\pi iq}$, we know that by \Cref{RootC}
\begin{eqnarray}
\nonumber
\mathcal{V}_{Log(\chi)} = 
 \begin{cases} 
      \mathcal{O}(0) & \text{when } q=0 \\
      \mathcal{O}(-1) & \text{when } q \not= 0.
 \end{cases}
\end{eqnarray}

$\underline{\text{Inductive Step}}$. Consider a representation of dimension $n$, $\rho: \mathbb{Z} \rightarrow \text{GL}_{n}(\mathbb{C})$ with $\rho$ reducible but indecomposable. 
It follows that $\rho$ is determined by $1 \mapsto M$, where $M$ is a Jordan block with eigenvalue $\lambda = re^{2\pi iq}$. 
Without loss of generality we may write $0 \rightarrow \chi_{q} \rightarrow \rho \rightarrow \rho' \rightarrow 0$ which by \Cref{Equiv-cat} and \Cref{Extension} we obtain $0 \rightarrow \mathcal{V}_{Log(\chi_{q})} \rightarrow \mathcal{V}_{Log(\rho)} \rightarrow \mathcal{V}_{Log(\rho')} \rightarrow 0$.
Moreover, $\mathcal{V}_{Log(\chi_{q})} \cong \mathcal{O}(\xi_{q})$ by the base case and $\mathcal{V}_{Log(\rho')} \cong \mathcal{O}(\xi_{q})^{\oplus(n-1)}$ by the inductive assumption. 

Hence, we have $0 \rightarrow \mathcal{O}(\xi_{q}) \rightarrow \mathcal{V}_{Log(\rho)} \rightarrow \mathcal{O}(\xi_{q})^{\oplus (n-1)} \rightarrow 0$.
Now, note that 
\begin{eqnarray}
\nonumber H^{1}(\mathbb{P}^{1}, (\mathcal{O}(\xi_{q})^{\oplus(n-1)})^{*} \otimes \mathcal{O}(\xi_{q})) && = H^{1}(\mathbb{P}^{1}, \mathcal{O}(-\xi_{q})^{\oplus(n-1)} \otimes \mathcal{O}(\xi_{q})) \\
\nonumber &&= H^{1}(\mathbb{P}^{1}, \mathcal{O}(0)^{\oplus(n-1)}).  
\end{eqnarray}

\noindent Further, by Serre-Duality we have
\begin{eqnarray}
\nonumber H^{1}(\mathbb{P}^{1}, \mathcal{O}(0)^{\oplus(n-1)}) &&= H^{0}(\mathbb{P}^{1}, (\mathcal{O}(0)^{\oplus(n-1)})^{*} \otimes \Omega_{\mathbb{P}^{1}}^{1})^{*} \\
\nonumber &&= H^{0}(\mathbb{P}^{1}, \mathcal{O}(-2)^{\oplus(n-1)})^{*}.
\end{eqnarray}

\noindent Of course dim $H^{0}(\mathbb{P}^{1}, \mathcal{O}(-2)^{\oplus(n-1)})^{*} = 0$ and so $0 \rightarrow \mathcal{O}(\xi_{q}) \rightarrow \mathcal{V}_{Log(\rho)} \rightarrow \mathcal{O}(\xi_{q})^{\oplus (n-1)} \rightarrow 0$ splits. 
Thus, we have the desired result which is 
\begin{eqnarray}
\nonumber
\mathcal{V}_{Log(\rho)} \cong 
 \begin{cases} 
      \mathcal{O}(0)^{\oplus n} & \text{when } q=0 \\
      \mathcal{O}(-1)^{\oplus n} & \text{when } q \not= 0.
 \end{cases}
\end{eqnarray}

For part $(2)$, let $\rho: \mathbb{Z} \rightarrow \text{GL}_{n}(\mathbb{C})$ be a monodromy representation of $\textbf{Y}_{(0,\infty)}$ with $(\mathcal{V}_{Log(\rho)}, \nabla_{Log(\rho)})$ the associated extended logarithmic connection where $\rho \cong \bigoplus_{i = 1}^{u} \rho_{i}$ such that each $\rho_{i}$ is indecomposable. 
As a result of \Cref{Equiv-cat} and \Cref{Extension}, we may write $\mathcal{V}_{\text{Log}(\rho)} \cong \bigoplus_{i = 1}^{u} \mathcal{V}_{\text{Log}(\rho_{i})}$. 
Since each $\rho_{i}$ is assumed to be indecomposable then we apply part $(1)$ above to obtain that for each $i, \mathcal{V}_{\text{Log}(\rho_{i})} \cong \mathcal{O}(-1)^{\oplus \text{ dim}(\rho_{i})}$ or $\mathcal{V}_{\text{Log}(\rho_{i})} \cong \mathcal{O}(0)^{\oplus \text{ dim}(\rho_{i})}$ dependent on the argument of the eigenvalue associated with the image of the generator. 
Now, allow $\mathcal{E}$ to be the direct sum of all the subbundles $\mathcal{V}_{\text{Log}(\rho_{j})}$ whose root is $-1$, then $\mathcal{E} = \bigoplus_{j = 1}^{\alpha} \mathcal{O}(-1)^{\oplus \text{ dim}(\rho_{j})} = \mathcal{O}(-1)^{\oplus \sum_{j = 1}^{\alpha} \text{ dim}(\rho_{j})} = \mathcal{O}(-1)^{\oplus -\text{c}_{1}(\mathcal{V}_{\text{Log}(\rho)})}$ where the last equality follows from \Cref{ChernR} combined with \Cref{RootC}. 
The remaining subbundles all have root $0$ and the dimension of their direct sum must be  $n + \text{c}_{1}(\mathcal{V}_{\text{Log}(\rho)})$. 
$\Box$

\begin{Rem} \label{IrreducibleZ}
   Recall that $\mathbb{Z}$ does not admit irreducible representations of dimension greater than one. 
   As a consequence, \Cref{RootC} along with \Cref{RootGen} close the case on monodromy representations of $\textbf{Y}_{(0,\infty)}$.
\end{Rem}

\section{Double Punctured Complex Line}

The goal of this section is to compute the roots arising from any two-dimensional reducible representation of $\textbf{Y}_{(0,1,\infty)}$.

\begin{Prop} \label{RootThreeC}
    Let $\chi: \mathbb{Z}*\mathbb{Z} \rightarrow \mathbb{C}^{*}$ be a monodromy representation of $\textbf{Y}_{(0, 1, \infty)}$ given by $\gamma_{0} \mapsto r_{0}e^{2\pi iq_{0}}$ and $\gamma_{1} \mapsto r_{1}e^{2\pi iq_{1}}$, with $(\mathcal{V}_{Log(\chi)}, \nabla_{Log(\chi)})$ the associated extended logarithmic connection. Then
    
    \begin{eqnarray}
    \nonumber
    \mathcal{V}_{Log(\chi)} \cong 
    \begin{cases} 
      \mathcal{O}(0) & \text{when } q_{j}=0 \\
      \mathcal{O}(-1) & \text{when } \sum q_{j} \leq 1 \text{ and not both } q_{j} = 0 \\
      \mathcal{O}(-2) & \text{when } \sum q_{j} > 1
    \end{cases}
    \end{eqnarray}.
\end{Prop}

\proof Observe that the connection map has three poles, one at $0$, one at $1$, and one at $\infty$ with respective residues $\text{log}(r_{0}e^{2 \pi iq_{0}})$, 
$\text{log}(r_{1}e^{2 \pi iq_{1}})$, 
$\text{log}((r_{0}r_{1})^{-1}e^{-2\pi i (q_{0} + q_{1})})$. 
By Ohtsuki’s formula we have that $c_{1}(\mathcal{V}_{Log(\chi)}) = -\text{log} (r_{0}e^{2 \pi iq_{0}}) -\text{log} (r_{1}e^{2 \pi iq_{1}})-\text{log}(\frac{1}{r_{0}r_{1}}e^{-2 \pi i(q_{0} + q_{1})})$, consequently, choosing the branch of log defined by $0 \leq q_{j} < 1$ we have the following set of cases. 

Case $1$. If $q_{j} = 0$ it is clear that $c_{1}(\mathcal{V}_{Log(\chi)}) = 0$.

Case $2$. Suppose that $\sum q_{j} \leq 1$ and either $q_{0} \not= 0$ or $q_{1} \not= 0$. Since by assumption $\sum q_{j} \leq 1$ and not both $q_{j}$ are zero, then $-\sum q_{j} \not\in [0,1)$ and so we must use $-\sum q_{j} + 1$, thus, $c_{1}(\mathcal{V}_{Log(\chi)}) = - q_{0} - q_{1} -(-q_{0} - q_{1} + 1) = -1$.

Case $3$. Suppose that $\sum q_{j} > 1$. Then, since $\sum q_{j} > 1$ it follows that $-\sum q_{j} \in (-2,-1)$, so we must use $-\sum q_{j} + 2 \in (0,1)$. Thus, $c_{1}(\mathcal{V}_{Log(\chi)}) = - q_{0} - q_{1} -(-q_{0} - q_{1} + 2) = -2$.

Observe that the result follows from \Cref{first-chern-class}.
$\Box$

\begin{Rem}
    Before introducing the main result of this section, we first prove the following lemma regarding the first Chern class of associated extended logarithmic connections that arise from monodromy representations with some specific constraints. 
\end{Rem}

\begin{Lm}\label{ChernThreeTwo}
    Suppose that $0 \rightarrow \rho' \rightarrow \rho \rightarrow \rho'' \rightarrow 0$ is a short exact sequence of monodromy representations of $\textbf{Y}_{(0, 1, \infty)}$ with $\rho$ two-dimensional. 
    Allow $\rho$ to be defined by $\gamma_{0} \mapsto M_{0}$ and $\gamma_{1} \mapsto M_{1}$ where $\lambda_{0} = r_{\lambda,0}e^{2\pi iq_{\lambda, 0}} \text{ and } \mu_{0} =r_{\mu,0}e^{2\pi iq_{\mu,0}}$ are the eigenvalues of $M_{0}$, while $\lambda_{1} = r_{\lambda,1}e^{2\pi iq_{\lambda, 1}} \text{ and } \mu_{1} =r_{\mu,1}e^{2\pi iq_{\mu,1}}$ are the eigenvalues of $M_{1}$ with all $q$ such that $0 \leq q < 1$. 
    Further, suppose that the $\lambda$s determine $\rho'$ and the $\mu$s determine $\rho''$.
    Then, if $\sum q_{\lambda, \iota} >1$ and $q_{\mu,j} = 0$, the first Chern class of the associated extended logarithmic connection $(\mathcal{V}_{Log(\rho)}, \nabla_{Log(\rho)})$ is $c_{1}(\mathcal{V}_{Log(\rho)}) = -2$.
\end{Lm}

\proof Note that by the classification of 2-dimensional representations of the free group with two generators found in \cite[Ch $2$, Theorem $4.2.1$]{IKSY2013}, there are three conjugacy classes of $\rho$ that fit the assumptions that $\rho$ is reducible, $\sum q_{\lambda, \iota} >1$, and $q_{\mu,j} = 0$ with $0 \leq q <1$. 
The three conjugacy classes are given by
\begin{eqnarray}
    \label{Reduced}
    \nonumber
    M_{0} \leftrightarrow 
    \begin{pmatrix}
        \lambda_{0} && 0 \\
        0 && \mu_{0} 
    \end{pmatrix}, \text{  }
    M_{1} \leftrightarrow 
    \begin{pmatrix}
        \lambda_{1} && 0 \\
        0 && \mu_{1} 
    \end{pmatrix}; \\
    \nonumber
    M_{0} \leftrightarrow 
    \begin{pmatrix}
        \lambda_{0} && 0 \\
        0 && \mu_{0} 
    \end{pmatrix}, \text{  }
    M_{1} \leftrightarrow 
    \begin{pmatrix}
        \lambda_{1} && 1 \\
        0 && \mu_{1} 
    \end{pmatrix}; \\
    \nonumber
    M_{0} \leftrightarrow 
    \begin{pmatrix}
        \mu_{0} && 0 \\
        0 && \lambda_{0} 
    \end{pmatrix}, \text{  }
    M_{1} \leftrightarrow 
    \begin{pmatrix}
        \mu_{1} && 1 \\
        0 && \lambda_{1} 
    \end{pmatrix}.
\end{eqnarray}

The connection map has three poles, one at $0$, one at $1$, and one at $\infty$ with respective residues $\text{log}(M_{0})$, 
$\text{log}(M_{1})$, 
$\text{log}((M_{0}M_{1})^{-1})$. 
By Ohtsuki’s formula, $c_{1}(\mathcal{V}_{Log(\rho)}) = -\text{Tr}(\text{log}(M_{0})) -\text{Tr}(\text{log}(M_{1})) -\text{Tr}(\text{log}((M_{0}M_{1})^{-1}))$. 
In the case where both $M_{0}$ and $M_{1}$ are diagonalizable it is clear that
\begin{eqnarray}
    \nonumber
    c_{1}(\mathcal{V}_{Log(\rho)}) = && - \text{log}(r_{\lambda,0}e^{2\pi iq_{\lambda, 0}}) - \text{log}(r_{\mu,0}e^{2\pi i(0)}) - \text{log}(r_{\lambda,1}e^{2\pi iq_{\lambda, 1}}) - \text{log}(r_{\mu,1}e^{2\pi i(0)}) \\
    \nonumber
    && - \text{log}\Big(\frac{e^{2\pi i(2-\sum q_{\lambda,\iota})}}{r_{\lambda,0}r_{\lambda,1}}\Big) - \text{log}\Big(\frac{e^{2\pi i(0)}}{r_{\mu,0}r_{\mu,1}}\Big)
\end{eqnarray}

\noindent where $2- \sum q_{\lambda,\iota}$ arises from the fact that $\sum q_{\lambda,\iota} > 1$ along with the choice of log as in the proof of \Cref{RootThreeC}. Moreover, note that all the $r$s cancel out. 
Thus, 
\begin{eqnarray}
  \nonumber
  c_{1}(\mathcal{V}_{Log(\rho)}) && = - \text{log}(\lambda_{0}) - \text{log}(\mu_{0}) - \text{log}(\lambda_{1}) - \text{log}(\mu_{1}) - \text{log}((\lambda_{0}\lambda_{-1})^{-1}) - \text{log}((\mu_{0}\mu_{1})^{-1}) \\
  \nonumber
  &&= - q_{\lambda,0} - q_{\lambda,1} - (-q_{\lambda,0} - q_{\lambda,1} + 2) \\
  \nonumber
  &&= -2.  
\end{eqnarray}

For the second case, again by Ohtsuki’s formula, $c_{1}(\mathcal{V}_{Log(\rho)}) = -\text{Tr}(\text{log}(M_{0})) -\text{Tr}(\text{log}(M_{1})) -\text{Tr}(\text{log}((M_{0}M_{1})^{-1}))$. 
The eigenvalues of log($M_{0}$) are log($\lambda_{0}$) and log($\mu_{0}$). 
By assumption $\sum q_{\lambda,\iota} > 1$, $q_{\mu,j} = 0$, and by the choice of log $0 \leq q < 1$ for all $q$. As a result $\lambda_{1} \not= \mu_{1}$ and hence there exists a basis $\Upsilon'$ where $M_{1}$ is diagonal. 
In the basis $\Upsilon'$, the trace of log$(M_{1})$ is $\text{log}(\lambda_{1}) + \text{log}(\mu_{1})$.
With regards to $(M_{0}M_{1})^{-1}$, we know that $M_{0}M_{1}$ is upper triangular with diagonal entries $\lambda_{0}\lambda_{1}$ and $\mu_{0}\mu_{1}$, further, as a consequence of  $\sum q_{\lambda, \iota} >1$ and $q_{\mu,j} = 0$, with $0 \leq q <1$, it follows that $\lambda_{0}\lambda_{1} \not= \mu_{0}\mu_{1}$.
Hence, there exists a basis, $\Upsilon''$, where $M_{0}M_{1}$ is diagonal. In the basis $\Upsilon''$, the trace of log($(M_{0}M_{1})^{-1}$) is $\text{log}((\lambda_{0}\lambda_{1})^{-1}) + \text{log}((\mu_{0}\mu_{1})^{-1})$.
Thus, this case collapses back into case 1. 
The third case follows the exact same procedure as the second case.
$\Box$

\begin{Thm} \label{NeedDer}
    Suppose that $0 \rightarrow \rho' \rightarrow \rho \rightarrow \rho'' \rightarrow 0$ is a short exact sequence of monodromy representations of $\textbf{Y}_{(0, 1, \infty)}$ with $\rho$ two-dimensional and $\rho', \rho''$ one-dimensional. 
    With the principal branch as a choice of log, let $0 \rightarrow \mathcal{O}(\xi') \rightarrow \mathcal{V}_{Log(\rho)} \rightarrow \mathcal{O}(\xi'') \rightarrow 0$ be the short exact sequence of the associated extended locally free sheaves with logarithmic connection maps.
    \begin{enumerate}
        \item If $\xi',\xi'' \not= -2,0$ respectively, then $0 \rightarrow \mathcal{O}(\xi') \rightarrow \mathcal{V}_{Log(\rho)} \rightarrow \mathcal{O}(\xi'') \rightarrow 0$ splits.
        \item If $\xi', \xi'' = -2,0$ respectively, then $\mathcal{V}_{Log(\rho)} \cong \mathcal{O}(-1) \oplus \mathcal{O}(-1)$ or $\mathcal{V}_{Log(\rho)} \cong \mathcal{O}(0) \oplus \mathcal{O}(-2)$.
    \end{enumerate} 
\end{Thm}

\proof Let $\rho: \mathbb{Z}*\mathbb{Z} \rightarrow \text{GL}_{2}(\mathbb{C})$ be a reducible 2-dimensional representation defined by $\gamma_{0} \mapsto M_{0}$ and $\gamma_{1} \mapsto M_{1}$, where $\lambda_{0}, \mu_{0}$ and $\lambda_{1}, \mu_{1}$ are the eigenvalues of $M_{0}$ and $M_{1}$ respectively. Further, let $\lambda_{0} = r_{\lambda,0}e^{2\pi i q_{\lambda,0}}, \mu_{0} = r_{\mu,0}e^{2\pi i q_{\mu,0}}, \lambda_{1} = r_{\lambda, 1} e^{2\pi i q_{\lambda,1}},$ and $\mu_{1} = r_{\mu,1}e^{2\pi i q_{\mu,1}}$. 
As $\rho$ is reducible, then by \cite[Ch $2$, Theorem $4.2.1$]{IKSY2013} there exists a basis $\hat{\Upsilon}$ in which
\begin{eqnarray}
    \nonumber
    M_{0} \leftrightarrow 
    \begin{pmatrix}
        \lambda_{0} && \eta_{0} \\
        0 && \mu_{0} 
    \end{pmatrix} \text{ and  }
    M_{1} \leftrightarrow 
    \begin{pmatrix}
        \lambda_{1} && \eta_{1} \\
        0 && \mu_{1} 
    \end{pmatrix}
\end{eqnarray}

\noindent and moreover, we may write $0 \rightarrow \rho' \rightarrow \rho \rightarrow \rho'' \rightarrow 0$, where $\rho'$ is determined by the $\lambda$s and $\rho''$ is determined by the $\mu$s.

For $(1)$, using cohomology, 
\begin{eqnarray}
\nonumber H^{1}(\mathbb{P}^{1}, (\mathcal{O}(\xi'')^{*} \otimes \mathcal{O}(\xi')) && = H^{1}(\mathbb{P}^{1}, \mathcal{O}(-\xi'') \otimes \mathcal{O}(\xi')) \\
\nonumber &&= H^{1}(\mathbb{P}^{1}, \mathcal{O}(-\xi''+\xi')).  
\end{eqnarray}

\noindent Further by Serre-Duality we have
\begin{eqnarray}
\nonumber H^{1}(\mathbb{P}^{1}, \mathcal{O}(-\xi''+\xi')) &&= H^{0}(\mathbb{P}^{1}, (\mathcal{O}(-\xi''+\xi')^{*} \otimes \Omega_{\mathbb{P}^{1}}^{1}))^{*} \\
\nonumber &&= H^{0}(\mathbb{P}^{1}, \mathcal{O}(\xi''-\xi'-2))^{*}.
\end{eqnarray}

\noindent Hence, if $\xi'' - \xi' < 2$ then $\mathcal{V}_{Log(\rho)} = \mathcal{O}(\xi') \oplus \mathcal{O}(\xi'')$. 
At this point observe that by \Cref{RootThreeC} we know that 
\begin{eqnarray}
    \nonumber
    \xi' = 
    \begin{cases} 
      0 & \text{when } q_{\lambda,i}=0 \\
      -1 & \text{when } \sum q_{\lambda,i} \leq 1 \text{ and not both } q_{\lambda,i} = 0 \\
      -2 & \text{when } \sum q_{\lambda,i} > 1
    \end{cases}
\end{eqnarray}

\noindent and similarly, 
\begin{eqnarray}
    \nonumber
    \xi'' = 
    \begin{cases} 
      0 & \text{when } q_{\mu,i}=0 \\
      -1 & \text{when } \sum q_{\mu,i} \leq 1 \text{ and not both } q_{\mu,i} = 0 \\
      -2 & \text{when } \sum q_{\mu,i} > 1.
    \end{cases}
\end{eqnarray}

\noindent It is clear that the only case that does not satisfy the inequality $\xi'' - \xi' < 2$, is when $\xi' = -2$ and $\xi'' = 0$.

For part $(2)$, the event that $\xi' = -2$ and $\xi'' = 0$ gives rise to the following short exact sequence $0 \rightarrow \mathcal{O}(-2) \rightarrow \mathcal{V}_{Log(\rho)} \rightarrow \mathcal{O}(0) \rightarrow 0$. 
Looking at the long exact sequence in cohomology 
\begin{eqnarray}
    \nonumber
    && 0 \rightarrow H^{0}(\mathbb{P}^{1}, \mathcal{O}(-2)) \rightarrow H^{0}(\mathbb{P}^{1}, \mathcal{V}_{Log(\rho)}) \rightarrow H^{0}(\mathbb{P}^{1}, \mathcal{O}(0)) \xrightarrow{\delta} H^{1}(\mathbb{P}^{1}, \mathcal{O}(-2)) \\ 
    \nonumber
    &&\rightarrow H^{1}(\mathbb{P}^{1}, \mathcal{V}_{Log(\rho)}) \rightarrow H^{1}(\mathbb{P}^{1}, \mathcal{O}(0)) \cdots
\end{eqnarray}

\noindent we note that $H^{0}(\mathbb{P}^{1}, \mathcal{O}(-2))$ and $H^{1}(\mathbb{P}^{1}, \mathcal{O}(0))$ are zero-dimensional complex vector spaces while $H^{1}(\mathbb{P}^{1}, \mathcal{O}(-2))$ is 1-dimensional. 
Hence, $\delta$ is either an isomorphism or the $0$ map. 

If $\delta$ is an isomorphism, then dim($H^{0}(\mathbb{P}^{1},\mathcal{V}_{Log(\rho)})) = 0$, meaning that $\mathcal{V}_{Log(\rho)} \cong \mathcal{O}(\xi_{1}) \oplus \mathcal{O}(\xi_{2})$ with $\xi_{1}, \xi_{2} \in \mathbb{Z}^{-}$. 
If $\delta$ is the $0$ map, then $H^{0}(\mathbb{P}^{1}, \mathcal{V}_{Log(\rho)}) \cong \mathbb{C}$, meaning that $\mathcal{V}_{Log(\rho)} \cong \mathcal{O}(0) \oplus \mathcal{O}(\xi)$ with $\xi \in \mathbb{Z}^{-}$. 

However, given that we are working in the scenario where $\sum q_{\lambda,i} > 1$ and $q_{\mu,i} = 0$, then by \Cref{ChernThreeTwo}, $c_{1}(\mathcal{V}_{Log(\rho)}) = -2$. Thus, if $\delta$ is an isomorphism then $\mathcal{V}_{Log(\rho)} \cong \mathcal{O}(-1) \oplus \mathcal{O}(-1)$, whereas if $\delta$ is the $0$ map, then $\mathcal{V}_{Log(\rho)} \cong \mathcal{O}(0) \oplus \mathcal{O}(-2)$.
$\Box$

\section{The Monodromy Derivative}

Examining the case of $\textbf{Y}_{(0, \infty)}$, we recall that $\mathbb{Z}$ does not admit irreducible representations of dimension greater than one, as was highlighted by \Cref{IrreducibleZ}. This meant that on every occasion we considered representations $\rho$ with dim$(\rho) > 1$ for which we were able to construct a short exact sequence $0 \rightarrow \rho' \rightarrow \rho \rightarrow \rho'' \rightarrow 0$. As a result of this short exact sequence we were allowed to apply sheaf cohomology on the resulting short exact sequence of associated extended logarithmic connections. 

When studying $\textbf{Y}_{(0, 1, \infty)}$, its fundamental group is $\mathbb{Z} * \mathbb{Z}$, which admits irreducible representations of dimension greater than $1$, particularly of dimension $2$ as seen in \cite[Ch $2$, Theorem $4.2.1$]{IKSY2013}. Thus, the sheaf cohomology approach is not a fruitful approach when $\rho$ is irreducible of dimension greater than one. 

With regard to the case when $\rho$ is irreducible, the idea is to construct a Serre-like derivative arising from $\textbf{Y}_{(0,1,\infty)}$. 

\subsection{Auxiliary Connection}

The construction of the Serre-like derivative will be two fold carried in this subsection and the next: concluding with \Cref{monodromy-derivative}. 
The purpose of this subsection is to develop what we will call an \textit{auxiliary connection} arising from $\textbf{Y}_{(0,1,\infty)}$.

Take the character representation $\chi: \pi_{1}(\mathbf{Y}_{(0,1,\infty)}, y) \rightarrow \mathbb{C}^{*}$ defined by mapping both generators to $-1$, applying \Cref{RootThreeC}, the resulting associated extended logarithmic connection is $(\mathcal{V}_{Log(\chi)}, \nabla_{Log(\chi)})$ over $\mathbb{P}^{1}$ with $\mathcal{V}_{Log(\chi)} \cong \mathcal{O}(-1)$.
Moreover, we dualize and obtain $(\mathcal{O}(1), \nabla_{\text{Log}(\chi)}^{*})$.
At this point we make a simple change in notation and refer to $\nabla_{\text{Log}(\chi)}^{*}$ simply as $\nabla$.
Now globally,
\begin{eqnarray}
    \nabla: \mathcal{O}(1) \rightarrow \mathcal{O}(1) \otimes_{\mathcal{O}} \Omega_{\mathbb{P}^{1}}^{1}([0]+[1]+[\infty])
\end{eqnarray}
\noindent where  
\begin{eqnarray}
    \nonumber
    \Omega_{\mathbb{P}^{1}}^{1}([0]+[1]+[\infty]) && = \Omega_{\mathbb{P}^{1}}^{1} \otimes_{\mathcal{O}} \mathcal{O}([0]+[1]+[\infty]) \\ 
    \nonumber
    && \cong \mathcal{O}(-2) \otimes_{\mathcal{O}} \mathcal{O}(3) \\
    \nonumber
    && \cong \mathcal{O}(1)
\end{eqnarray}

\noindent and so
\begin{eqnarray}
    \nabla: \mathcal{O}(1) \rightarrow \mathcal{O}(2).
\end{eqnarray}

Next, we take $\nabla^{\otimes 2}$ to be the connection map on $\mathcal{O}(1) \otimes \mathcal{O}(1)$ such that for all sections $s,t \in O(1)$ the following is satisfied: $(\nabla \otimes \nabla)(s \otimes t) = \nabla(s) \otimes t + s \otimes \nabla(t)$. 
Observe that 
\begin{eqnarray}
    \nabla^{\otimes 2}: \mathcal{O}(2) \rightarrow \mathcal{O}(2) \otimes \Omega_{\mathbb{P}^{1}}^{1}([0]+[1]+[\infty]),
\end{eqnarray}
 
\noindent and so, inductively,  we define $\nabla_{\xi} := \nabla^{\otimes \xi} : \mathcal{O}(\xi) \rightarrow \mathcal{O}(\xi) \otimes_{\mathcal{O}} \Omega^{1}([0]+[1]+[\infty])$. 
Hence, using the fact that $ \Omega_{\mathbb{P}^{1}}^{1}([0]+[1]+[\infty]) \cong \mathcal{O}(1)$, then
\begin{eqnarray}
    \nabla_{\xi}: \mathcal{O}(\xi) \rightarrow \mathcal{O}(\xi+1).
\end{eqnarray}

\begin{Defn} \label{little-delta}
    In the construction above we refer to $\nabla_{\xi}$ as the \textit{auxiliary connection map of weight} $\xi$ \textit{arising from} $\textbf{Y}_{(0,1,\infty)}$.
\end{Defn}

\begin{Rem} \label{general-little-delta}
     The construction of the auxiliary connection map can be generalized as arising from the $m$th punctured projective line. 
     However, observe that as the number of removed points increases so does the degree of the auxiliary connection arising from the space. To see this explicitly, let $p_{i}$ be points of $\mathbb{P}^{1}$, then 
     \begin{eqnarray}
         \nonumber
         \Omega^{1}_{\mathbb{P}^{1}} \big(\sum_{i=1}^{m}[p_{i}] \big) && = \Omega_{\mathbb{P}^{1}}^{1} \otimes_{\mathcal{O}} \mathcal{O} \big(\sum_{i=1}^{m}[p_{i}] \big) \\
         \nonumber
         && \cong \mathcal{O}(-2) \otimes_{\mathcal{O}} \mathcal{O}(m) \\
         \nonumber
         && \cong \mathcal{O}(m-2)
     \end{eqnarray}

     \noindent and so arising from $\mathbb{P}^{1}$ minus $m$ points 
     \begin{eqnarray}
         \nabla: \mathcal{O}(1) \rightarrow \mathcal{O}(m-1).
     \end{eqnarray}

     Hence, which is why when the auxiliary connection map arises from $\textbf{Y}_{(0,1,\infty)}$ the co-domain is $\mathcal{O}(2)$ and whereas if it arises from $\textbf{Y}_{(0,\infty)}$ the co-domain would be $\mathcal{O}(1)$.
\end{Rem}

\begin{Ex} \label{little-delta-example}
    Consider $\textbf{Y}_{(0, \infty)}$ whose fundamental group is $\mathbb{Z}$. 
    Let $\chi$ be a character representation that maps the generator to $-1$, then after dualizing, the associated extended vector bundle at hand is $\mathcal{O}(1)$.
    It follows the that the auxiliary connection arising from $\textbf{Y}_{(0,\infty)}$ is a map such that 
    \begin{eqnarray}
        \nabla: \mathcal{O}(1) \rightarrow \mathcal{O}(1) \otimes_{\mathcal{O}} \Omega^{1}([0]+[\infty]).
    \end{eqnarray}

    \noindent Observe that $\Omega^{1}([0] + [\infty]) \overset{\phi}{\cong} \mathcal{O}(0)$ by $\phi(\frac{dz}{z}) = 1$, which means that 
    \begin{eqnarray}
        \nabla: \mathcal{O}(1) \rightarrow \mathcal{O}(1).
    \end{eqnarray}
    
    Over global sections, $\nabla$ gives rise to a degree preserving map $\delta: \mathbb{C}[x,y] \rightarrow \mathbb{C}[x,y]$, where $x$ has a zero at $0$, and $y$ has a zero at $\infty$. 
    Now, since $\delta$ is a degree preserving map, we must have $\delta x = a_{1}x + a_{2}y$ and $\delta y = b_{1}x + b_{2}y$ with $a_{j}, b_{i} \in \mathbb{C}$. 

    Let $U_{0} := \mathbb{P}^{1} - \{0\}$ and $U_{\infty} := \mathbb{P}^{1} - \{\infty\}$. Restricting over $U_{\infty}$ gives, Res$(\delta x) = \text{Res}(a_{1}x + a_{2}y)$, and hence, $\delta z = a_{1}z + a_{2}$ as $\delta$ commutes with restriction, where $z = \frac{x}{y}$. 
    Moreover, as locally over $\textbf{Y}_{(0,\infty)}$, $\nabla = d - \frac{1}{2}\frac{dz}{z}$, then
    \begin{eqnarray}
        \nonumber
        \delta z &&= dz - \frac{dz}{2z}(z) \\
        \nonumber
        && = z\frac{dz}{z} - \frac{z}{2}\frac{dz}{z} \\
        \nonumber
        && = (\frac{1}{2}z)\frac{dz}{z}
    \end{eqnarray}
    however, after applying $\phi$ we are left with $\delta z = \frac{1}{2}z$.
    At this point we compare coefficients, $a_{1} = \frac{1}{2}$ and $a_{2} = 0$. 
    Similarly for $U_{0}$, we arrive at $b_{1} = 0$ and $b_{2} = -\frac{1}{2}$. Thus, $\delta x = \frac{1}{2}x$ and $\delta y = - \frac{1}{2} y$.
\end{Ex}

\subsection{Monodromy Derivative}

We now return to the case at hand, that is, considering irreducible monodromy representations of $\textbf{Y}_{(0,1,\infty)}$ of dimension $2$. Justifying the naming in \Cref{little-delta} we have the following: given a monodromy representation $\rho$, the associated extended logarithmic connection $(\mathcal{V}_{Log(\rho)}, \nabla_{Log(\rho)})$ can be seen as a triple, $(\mathcal{V}_{Log(\rho)}, \nabla_{Log(\rho)}, \nabla)$, where  $\nabla$ is the auxiliary connection map arising from $\textbf{Y}_{(0,1,\infty)}$.

Thus, keeping in mind that $\Omega_{\mathbb{P}^{1}}^{1}([0]+[1]+[\infty]) \cong \mathcal{O}(1)$, we now have the extra tool 
\begin{eqnarray}
    \nabla_{Log(\rho)} \otimes \nabla: \mathcal{V}_{Log(\rho)} \otimes_{\mathcal{O}} \mathcal{O}(1) \rightarrow \mathcal{V}_{Log(\rho)} \otimes_{\mathcal{O}} \mathcal{O}(2) 
\end{eqnarray}

\noindent and more generally, considering an auxiliary connection map of arbitrary weight 
\begin{eqnarray}
    \nabla_{Log(\rho)} \otimes \nabla_{\xi}: \mathcal{V}_{Log(\rho)} \otimes_{\mathcal{O}} \mathcal{O}(\xi) \rightarrow \mathcal{V}_{Log(\rho)} \otimes_{\mathcal{O}} \mathcal{O}(\xi+1). 
\end{eqnarray}

\begin{Defn} \label{twisted-module}
    Let $\mathcal{N}_{\xi}(\rho) := H^{0}(\mathbb{P}^{1}, \mathcal{V}_{\text{Log}(\rho)} \otimes_{\mathcal{O}} \mathcal{O}(\xi))$, and furthermore, define $\mathcal{N}(\rho) := \bigoplus_{\xi \in \mathbb{Z}} \mathcal{N}_{\xi}(\rho)$.
    We refer to $\mathcal{N}(\rho)$ as the \textit{twisted module of $\rho$}.
\end{Defn}

Since $\mathcal{V}_{\text{Log}(\rho)}$ is a direct sum of twisting sheaves it is easy to see that there exists a $\xi_{\text{min}}$ such that for any $\xi_{j} < \xi_{\text{min}}$ we have  $\mathcal{N}_{\xi_{j}}(\rho) = 0$.
We refer to $\xi_{\text{min}}$ as the \textit{minimal weight} of $\mathcal{N}(\rho)$.

Moreover, $\mathcal{N}(\rho)$ is a $\mathbb{Z}$-graded module of global sections over the ring $R := \mathcal{N}(\textbf{1}) = \mathbb{C}[x,y]$ where, of course, $\mathbb{P}^{1} = \text{Proj } \mathbb{C}[x,y]$ and $\textbf{1}$ is the trivial one-dimensional representation. 
By the Birkhoff-Grothendieck theorem along with the fact that $H^{0}$ and direct sums commute, this module is free of rank equal to the dimension of the representation; indeed, 
\begin{eqnarray}
    \mathcal{N}(\rho) &&= \bigoplus_{\xi \in \mathbb{Z}} H^{0}(\mathbb{P}^{1}, \mathcal{V}_{\text{Log}(\rho)} \otimes_{\mathcal{O}} \mathcal{O}(\xi)) \\
    && \cong \bigoplus_{\xi \in \mathbb{Z}} H^{0}(\mathbb{P}^{1}, \bigoplus_{i=1}^{\text{dim}(\rho)} \mathcal{O}(\xi_{i} + \xi)) \\
    && \cong \bigoplus_{\xi \in \mathbb{Z}} \bigoplus_{i=1}^{\text{dim}(\rho)} H^{0}(\mathbb{P}^{1},  \mathcal{O}(\xi_{i} + \xi)) \\
    && \cong \bigoplus\limits_{i=1}^{\text{dim}(\rho)} R[-\xi_{i}]
\end{eqnarray}
thus, any choice of homogeneous generators gives an isomorphism
\begin{eqnarray}
    \mathcal{N}(\rho) \cong \bigoplus\limits_{i=1}^{\text{dim}(\rho)} R[-\xi_{i}]
\end{eqnarray}

\noindent where by $R[a]$ we denote the rank one graded module over $R$ obtained by shifting the grading by $a$. 
On global sections, note that the new connection map $\nabla_{Log(\rho)} \otimes \nabla_{\xi}$ gives us a graded derivation of degree one. 
Indeed, let us denote the connection map $\nabla_{Log(\rho)} \otimes \nabla_{\xi}$ as $D_{\xi}$ when applied to global sections, so that
\begin{eqnarray}
    D_{\xi}: \mathcal{N}_{\xi}(\rho) \rightarrow \mathcal{N}_{\xi+1}(\rho)
\end{eqnarray}

\noindent where of course $D_{\xi}$ satisfies the Leibniz rule as $\nabla_{Log(\rho)} \otimes \nabla_{\xi}$ is a connection map and so 
\begin{eqnarray}
    D := \bigoplus_{\xi \in \mathbb{Z}} D_{\xi}
\end{eqnarray}

\noindent is a graded derivation of degree one on $\mathcal{N}(\rho)$.

\begin{Defn} \label{monodromy-derivative}
    We refer to the the graded derivation $D$ as the \textit{monodromy derivative of} $\rho$.
\end{Defn}

\begin{Rem}
    Analogous to \Cref{general-little-delta}, it is not hard to see how to generalize the monodromy derivative not only to any dimension of a monodromy representation of $\textbf{Y}_{(0,1,\infty)}$, but to any monodromy representation of an $m$th punctured projective line. We provide a rough sketch of how to generalize.
    Let $\textbf{Y}_{(m\text{th})}$ be the $m$th punctured projective line. 
    Given a monodromy representation of $\textbf{Y}_{(m\text{th})}$ we obtain the triplet $(\mathcal{V}_{\text{Log}(\rho)}, \nabla_{\text{Log}(\rho)}, \nabla)$ where $\nabla$ is the auxiliary connection map arising from $\textbf{Y}_{(m\text{th})}$.
    Next, we tensor the connection maps.
    However, as mentioned in \Cref{general-little-delta} the co-domain of the connection map $\nabla_{\text{Log}(\rho)} \otimes \nabla$ changes according to the number of removed points, hence, 
    \begin{eqnarray} 
        \nabla_{\text{Log}(\rho)} \otimes \nabla : \mathcal{V}_{Log(\rho)} \otimes_{\mathcal{O}} \mathcal{O}(1) \rightarrow \mathcal{V}_{Log(\rho)} \otimes_{\mathcal{O}} \mathcal{O}(m-1).
    \end{eqnarray}

    \noindent We can take global sections as before and create $\mathcal{N}(\rho)$. 
    Further, denoting the connection map $\nabla_{\text{Log}(\rho)} \otimes \nabla_{\xi}$ as $D_{\xi}$ when applied to global sections and defining $D := \bigoplus_{\xi \in \mathbb{Z}} D_{\xi}$ we have a graded derivation on $\mathcal{N}(\rho)$.
\end{Rem}

\begin{Rem}
    Now that we have constructed the monodromy derivative we will use it to prove the main theorem of this section, but first we direct our attention to the following two lemmas which will be used in the proof of the main theorem. 
\end{Rem}

\begin{Lm} \label{not-zero-derivative}
    Let $\xi_{\text{min}}$ be the minimal weight of the twisted module $\mathcal{N}(\rho)$ for some $n$-dimensional irreducible monodromy representation $\rho$ of $\textbf{Y}_{(0,1,\infty)}$ with $n \geq 2$. 
    Let $\xi_{\text{min}} \leq k$ and suppose that $f \in \mathcal{N}_{k}(\rho)$ with $f \not= 0$, then $D_{k}(f) \not= 0$.
\end{Lm}

\proof For the sake of contradiction let us assume that $D_{k}(f) = 0$. 

Restricting to $\textbf{Y}_{(0,1,\infty)}$, we see that $f|_{\textbf{Y}_{(0,1,\infty)}}$ is a global section of the associated complex bundle $\mathcal{Q}_{\rho}$ (see \Cref{Associated-Con}) because $\Big(\mathcal{V}_{\text{Log}(\rho)} \otimes \mathcal{O}(k)\Big)|_{\textbf{Y}_{(0,1,\infty)}} \cong \mathcal{Q}_{\rho}$.
Denote with $\mathcal{F} \subset \mathcal{Q}_{\rho}$ the subbundle generated by $f|_{\textbf{Y}_{(0,1,\infty)}}$. 
It is clear that the following diagram commutes  
\begin{center}
    \begin{tikzcd}
    \mathcal{V}_{\text{Log}(\rho)} \otimes \mathcal{O}(k) \arrow[r, "\nabla_{\text{Log}(\rho)} \otimes \nabla"] \arrow[d, "\varphi"]
    & \mathcal{V}_{\text{Log}(\rho)} \otimes \mathcal{O}(k+1) \arrow[d, "\varphi"] \\
    \mathcal{Q}_{\rho} \arrow[r, "(\nabla_{\text{Log}(\rho)} \otimes \nabla)|_{\textbf{Y}_{(0,1,\infty)}}"]
    & \mathcal{Q}_{\rho}  
    \end{tikzcd}
\end{center}

\noindent where $\varphi$ denotes the restriction map.

Since by assumption $D_{k}(f) = 0$, then after restricting $D_{k}|_{{\textbf{Y}_{(0,1,\infty)}}}(f|_{\textbf{Y}_{(0,1,\infty)}}) = 0$. 
Moreover, by \Cref{Free-sheaf}, $\mathcal{Q}_{\rho}$ is free over $\textbf{Y}_{(0,1,\infty)} = \text{Spec}(A)$ where $A = \mathbb{C}[x, x^{-1}, (x-1)^{-1}]$. 
Hence, by \Cref{Free-sheaf}, $\mathcal{F}$ can be viewed as a free module of rank one over $A$ whose generator is stable under $D|_{\textbf{Y}_{(0,1,\infty)}}$. 
Ergo, as $\mathcal{F}$ is stable under $D|_{\textbf{Y}_{(0,1,\infty)}} \cong \nabla_{\rho}$ then $\mathcal{F}$ inherits a connection map $\hat{\nabla}$ with $(\mathcal{F}, \hat{\nabla}) \subset (\mathcal{Q}_{\rho}, \nabla_{\rho})$ as connections.
Now, since there is an equivalence of categories between holomorphic connections on $\textbf{Y}_{(0,1,\infty)}$ and monodromy representations of $\textbf{Y}_{(0,1,\infty)}$, it follows that $\rho' \subset \rho$ where $\rho'$ corresponds to $(\mathcal{F}, \hat{\nabla})$. 
At this point notice that $(\mathcal{F}, \hat{\nabla})$ is nonempty which means that $\rho' \not= 0$, moreover, since $\mathcal{N}(\rho)$ is of rank $n$ it follows that $(\mathcal{F}, \hat{\nabla}) \subsetneqq (\mathcal{Q}_{\rho}, \nabla_{\rho})$ which implies that $\rho' \subsetneqq \rho$. 
Consequently, $\rho'$ must be a proper subrepresentation of $\rho$, but this is impossible since $\rho$ is irreducible. 
$\Box$

\begin{Lm} \label{main-theorem}
    Suppose that $\rho: \mathbb{Z}*\mathbb{Z} \rightarrow \text{GL}_{2}(\mathbb{C})$ is an irreducible two-dimensional monodromy representation of $\textbf{Y}_{(0,1,\infty)}$. 
    Denote the minimal weight of $\mathcal{N}(\rho)$ by $\xi_{\text{min}}$ and let $(\mathcal{V}_{\text{Log}(\rho)}, \nabla_{\text{Log}(\rho)})$ denote the associated extended logarithmic connection.
    Then either 
    \begin{eqnarray}
        \nonumber
        \mathcal{V}_{\text{Log}(\rho)} && \cong \mathcal{O}(-\xi_{\text{min}})^{\oplus 2} \text{ or} \\
        \nonumber
        \mathcal{V}_{\text{Log}(\rho)} && \cong \mathcal{O}(-\xi_{\text{min}}) \oplus \mathcal{O}(-\xi_{\text{min}}-1)
    \end{eqnarray}
\end{Lm}

\proof Since $\rho$ is two-dimensional $\mathcal{N}(\rho)$ has rank two, so let $F$ and $G$ be a pair of generators. 
Further, let $\xi_{\text{min}}, \xi_{1}$ be the weights of $F,G$ respectively, where by assumption $\xi_{\text{min}}$ is the minimal weight of $\mathcal{N}(\rho)$ and $\xi_{\text{min}} \leq \xi_{1}$.
Observe that by construction $D_{\xi_{\text{min}}}(F) \in \mathcal{N}_{\xi_{\text{min}}+1}(\rho)$ and so $D_{\xi_{\text{min}}}(F)$ has weight $\xi_{\text{min}} +1$. Moreover since $F$ and $G$ generate $\mathcal{N}(\rho)$ we must have $D_{\xi_{\text{min}}}(F) = \eta F + \tau G$ for some $\eta, \tau \in \mathbb{C}[x,y]$.
As a result of \Cref{not-zero-derivative} it follows that $D_{\xi_{\text{min}}}(F) \not= 0$.
At this point we remark that $\tau \not= 0$, as otherwise $F$ is stable under $D$ and an argument similar to that in the proof of \Cref{not-zero-derivative} will lead to a contradiction. 
Ergo, $\tau G$ must have weight $\xi_{\text{min}} + 1$.
Thus, $\xi_{\text{min}} + 1 = a + \xi_{1}$ implying that 
\begin{eqnarray}
    \nonumber
    \xi_{1} = 
    \begin{cases}
        \xi_{\text{min}} + 1 & \text{when } a=0 \\
        \xi_{\text{min}} & \text{when } a=1
    \end{cases}
\end{eqnarray}
 
\noindent since $0 \leq a$ with $\xi_{\text{min}} \leq \xi_{1}$.
$\Box$

\begin{Thm} \label{even-odd}
    Suppose that $\rho: \mathbb{Z}*\mathbb{Z} \rightarrow \text{GL}_{2}(\mathbb{C})$ is an irreducible two-dimensional monodromy representation of $\textbf{Y}_{(0,1,\infty)}$ and let $(\mathcal{V}_{\text{Log}(\rho)}, \nabla_{\text{Log}(\rho)})$ denote the associated extended logarithmic connection. 
    Allow $\zeta = c_{1}(\mathcal{V}_{\text{Log}(\rho)})$.
    Then 
    \begin{eqnarray}
        \nonumber
        \mathcal{V}_{\text{Log}(\rho)} \cong
        \begin{cases}
        \mathcal{O}(\frac{\zeta}{2})^{\oplus 2} & \text{when } \zeta \text{ is even} \\
        \mathcal{O}(\frac{\zeta - 1}{2}) \oplus \mathcal{O}(\frac{\zeta + 1}{2}) & \text{when } \zeta \text{ is odd}. 
    \end{cases}
    \end{eqnarray}
\end{Thm}

\proof Let $\zeta = 2\sigma$ for some $\sigma \in \mathbb{Z}$, and  where by \Cref{main-theorem}, for the sake of contradiction, we assume that $\mathcal{V}_{\text{Log}(\rho)} \cong \mathcal{O}(-\xi_{\text{min}}) \oplus \mathcal{O}(-\xi_{\text{min}}-1)$. 
Then, by \Cref{first-chern-class} it follows that $2\sigma = -2\xi_{\text{min}} - 1$ which is impossible as $\sigma, \xi_{\text{min}} \in \mathbb{Z}$. 
Hence, if $\zeta$ is even then $\mathcal{V}_{\text{Log}(\rho)} \cong \mathcal{O}(-\xi_{\text{min}})^{\oplus 2}$. 
Moreover, by \Cref{first-chern-class} the first Chern class is the sum of the roots, so $-\xi_{\text{min}} = \frac{\zeta}{2}$.
A similar argument shows the desired result for when $\zeta$ is odd. 
$\Box$

\begin{Rem}
    Witness that the theorem above does not make any assumptions regarding $\rho$ aside from the fact that it is an irreducible representation of dimension $2$.
    Thus, this theorem along with \Cref{NeedDer} close the case on two-dimensional monodromy representations of $\textbf{Y}_{(0,1,\infty)}$.
\end{Rem}

\section{Example from the Modular Group}

In this section we will explicitly compute the roots of an irreducible two-dimensional representation of $\mathbb{Z}*\mathbb{Z}$ that factors through $\text{SL}_{2}(\mathbb{Z})$ using the main result from the monodromy derivative. To remind the reader of such representations we closely follow \cite[\S $3$]{Mason2008}.

We introduce standard matrices 
\begin{eqnarray}
    \nonumber
    R =
    \begin{pmatrix}
        -1 && 1 \\
        -1 && 0 
    \end{pmatrix}, \text{  }
    S = 
    \begin{pmatrix}
        0 && -1 \\
        1 && 0 
    \end{pmatrix}, \text{  }
    T = 
    \begin{pmatrix}
        1 && 1 \\
        0 && 1 
    \end{pmatrix},
\end{eqnarray}

\noindent of which any two of them generate $\text{SL}_{2}(\mathbb{Z})$; indeed 
\begin{eqnarray}
    \nonumber 
    \Gamma := \text{SL}_{2}(\mathbb{Z}) = \langle R,S,T | RS = T, R^{3} = S^{4} = I \rangle.
\end{eqnarray}

\noindent Further, we assume that given $\rho: \Gamma \rightarrow \text{GL}_{2}(\mathbb{C})$ the matrix $\rho(T)$ is semisimple and $\rho$ is irreducible. For such representations, as \cite{Mason2008} shows, there exists a basis $\beta$ in which the images of $T,S$ and $R$ are not scalar matrices and
\begin{eqnarray}
    \nonumber 
    \rho(T) &&=
    \begin{pmatrix}
        \lambda && 0 \\
        0 && \lambda^{-1}\sigma 
    \end{pmatrix}, \text{  } \\
    \nonumber
    \rho(S) &&= 
    \begin{pmatrix}
        a && 1 \\
        -a^{2}-\sigma^{3} && -a 
    \end{pmatrix}, \text{  } \\
    \nonumber
    \rho(R) &&= \sigma^{3} 
    \begin{pmatrix}
        -a\lambda && -\lambda \\
        \sigma\lambda^{-1}(a^{2}+\sigma^{3}) && a\sigma\lambda^{-1} 
    \end{pmatrix}, 
\end{eqnarray}

\noindent where $\sigma$ is a sixth root of unity, $\sigma \not= \lambda^{2}$, $a = \frac{1}{\sigma(\lambda-\sigma\lambda^{-1})}$, and $a^{2} \not= -\sigma^{3}$.

In the basis $\beta$ above, we choose $\lambda = 1$, and $\sigma = -1$. 
    Observe then that the restrictions are satisfied and 
    \begin{eqnarray}
        \nonumber
        \rho(T) = 
        \begin{pmatrix}
            1 && 0 \\
            0 && -1
        \end{pmatrix}, \text { } 
        \rho(S) = 
        \begin{pmatrix}
            -\frac{1}{2} && 1 \\
            \frac{3}{4} && \frac{1}{2}
        \end{pmatrix}, \text{ }
        (\rho(T)\rho(S))^{-1} = 
        \begin{pmatrix}
            -\frac{1}{2} && -1 \\
            \frac{3}{4} && -\frac{1}{2}
        \end{pmatrix}.
    \end{eqnarray}
Now, as stated above, $\rho(T) \text{ and } \rho(S)$ give us an irreducible representation of $\Gamma$ and hence of $\mathbb{Z}*\mathbb{Z}$.
In other words, $\rho$ is an irreducible two-dimensional monodromy representation of $\textbf{Y}_{(0,1,\infty)}$. 
Thus, the representation above gives us the associated extended logarithmic connection $(\mathcal{V}_{Log(\rho)}, \nabla_{Log(\rho)})$ over $\mathbb{P}^{1}$.
The connection has three poles, one at $0$, one at $1$, and one at $\infty$ with respective residues $\text{log}(\rho(T))$, $\text{log}(\rho(S))$, $\text{log}((\rho(T)\rho(S))^{-1})$.

By Ohtsuki’s formula we have 
    \begin{eqnarray}
        \nonumber
        c_{1}(\mathcal{V}_{\text{Log}(\rho)}) = -\text{Tr}(\text{log}(\rho(T))) -\text{Tr}(\text{log}(\rho(S))) -\text{Tr}(\text{log}((\rho(T)\rho(S))^{-1}))
    \end{eqnarray}
\noindent and as a result of Equation  $(\ref{Jf})$ combined with the fact that the trace is invariant under conjugation we may use the Jordan form of each image.
Denote by $B_{J}$ the Jordan form of the matrix $B$.
    First, we have that $\rho(S) = H\rho(S)_{J}H^{-1}$ where 
    \begin{eqnarray}
        \rho(S)_{J}=
        \begin{pmatrix}
            e^{2\pi i(\frac{1}{2})} && 0 \\
            0 && e^{2 \pi i(0)} 
        \end{pmatrix} \text{with }
        H=
        \begin{pmatrix}
            -2 && \frac{2}{3} \\
            1 && 1 
        \end{pmatrix},
    \end{eqnarray}
    and further, we have $(\rho(T)\rho(S))^{-1} = K(\rho(T)\rho(S))^{-1}_{J}K^{-1}$ where
    \begin{eqnarray}
        (\rho(T)\rho(S))^{-1}_{J}=
        \begin{pmatrix}
            e^{2\pi i(\frac{2}{3})} && 0 \\
            0 && e^{2 \pi i(\frac{1}{3})} 
        \end{pmatrix} \text{with }
        K=
        \begin{pmatrix}
            -\frac{2i}{\sqrt{3}} && \frac{2i}{\sqrt{3}} \\
            1 && 1 
        \end{pmatrix}.
    \end{eqnarray}

    \noindent It follows then that 
    \begin{eqnarray}
        \nonumber
        \text{Tr}(\text{log}(\rho(T)_{J})) = \frac{1}{2}\\
        \nonumber
        \text{Tr}(\text{log}(\rho(S)_{J})) = \frac{1}{2}\\
        \nonumber 
        \text{Tr}(\text{log}((\rho(T)\rho(S))^{-1}_{J})) = 1
    \end{eqnarray}

    \noindent which implies that $c_{1}(\mathcal{V}_{Log(\rho)}) = -2$.
    Thus, by \Cref{even-odd} we must be in the case where
    \begin{eqnarray}
        \mathcal{V}_{Log(\rho)} \cong \mathcal{O}(-1)^{\oplus 2}.
    \end{eqnarray}

\bibliographystyle{siam}
\bibliography{References}

\end{document}